\documentclass[preprint]{elsarticle}

\usepackage{lineno,hyperref}
\usepackage{amsmath}
\usepackage{amsfonts}
\usepackage{graphicx}
\usepackage{subfig}
\usepackage[ruled,linesnumbered]{algorithm2e}

\usepackage{float}

\journal{Journal of Computational Physics}

\begin{document}

\begin{frontmatter}
\title{A FFT-based GMRES for fast solving of Poisson equation in concatenated geometry}




\author[1]{Zichao Jiang\corref{cor1}}
\cortext[cor1]{Corresponding author}
\ead{jiangzch8@mail.sysu.edu.cn (jzc_focal@hotmail.com)}
\author[1]{Jiacheng Lian}
\author[1]{Zhuolin Wang}
\address[1]{School of Aeronautics and Astronautics, Sun Yat-sen University, Shenzhen 518107, China}

\begin{abstract}
    Fast Fourier Transform (FFT)-based solvers for the Poisson equation are highly efficient, exhibiting $O(N\log N)$ computational complexity and excellent parallelism. However, their application is typically restricted to simple, regular geometries due to the separability requirement of the underlying discrete operators. This paper introduces a novel domain decomposition method that extends the applicability of FFT-based solvers to complex composite domains—geometries constructed from multiple sub-regions. The method transforms the global problem into a system of sub-problems coupled through Schur complements at the interfaces. A key challenge is that the Schur complement disrupts the matrix structure required for direct FFT-based inversion. To overcome this, we develop a FFT-based preconditioner to accelerate the Generalized Minimal Residual (GMRES) method for the interface system. The central innovation is a novel preconditioner based on the inverse of the block operator without the Schur complement, which can be applied efficiently using the FFT-based solver. The resulting preconditioned iteration retains an optimal complexity for each step. Numerical experiments on a cross-shaped domain demonstrate that the proposed solver achieves the expected second-order accuracy of the underlying finite difference scheme. Furthermore, it exhibits significantly improved computational performance compared to a classic sparse GMRES solver based on Eigen libeary. For a problem with $10^5$ grid points, our method achieves a speedup of over 40 times.
\end{abstract}


\end{frontmatter}

\section{Introduction}
\label{sec:intro}


This paper focuses on the Helmholtz-Poisson equation of the form
\begin{equation}
    \label{eq:Poisson}
        \nabla^2p+\kappa p=f,\boldsymbol{x}\in\Omega,
\end{equation}
where $\boldsymbol{x}=\left(x,y,z\right)$ represents the spatial coordinates, and $f$ is a given function on $\Omega$. The equation presented in \eqref{eq:Poisson} can be transformed into various other equations, including the Schrödinger equation. When the constant $\kappa=0$, the equation reduces to the Poisson equation, the solution of which holds significant importance for various applications, including incompressible flow simulations. Furthermore, the discrete scheme for this elliptic problem can be integrated into time-stepping schemes for parabolic and hyperbolic equations, such as the heat conduction and wave equations.

Over several decades, researchers have proposed numerous algorithms to enhance the computational efficiency of solving the equation \eqref{eq:Poisson}. Among these, the most representative are algorithms based on the FFT\cite{jodra2017solving, cooley1965algorithm}. Originally developed for solving the Poisson equation on rectangular domains\cite{skollermo1975fourier}, these methods were subsequently widely applied to various large-scale problems, including turbulence simulation\cite{xie2021low, he2024turbulent}. Their widespread adoption is largely attributed to their computational complexity of $O\left(N \log N\right)$, where $N$ denotes the total number of discrete grid points. Furthermore, these methods exhibit inherent parallelism, which makes them highly suitable for modern parallel computing architectures.

These methods typically employ discretizations based on the finite difference method (FDM) or the finite volume method (FVM)\cite{han2025convergent, coco2024high, liu2025sixth}. Initially, the global equations were transformed into independent sub-equations by expanding the solution in Fourier series, which were then solved using the Thomas algorithm. These approaches offer significant advantages in computational efficiency compared to Krylov subspace solvers and multigrid methods. However, these methods have a notable limitation. Their reliance on a spectral basis requires the discrete operator to be separable, a condition typically met only on structured, uniform grids. Consequently, these traditional methods makes the proposed methods challenging to handle problems with curved boundaries or complex geometries.

In recent years, to overcome the geometric adaptability limitations of traditional FFT algorithms, researchers have conducted extensive exploratory work. For instance, the FastRK3 solver \cite{aithal2020fast} employs a discretization scheme on general curvilinear coordinates, as proposed by Nikitin \cite{nikitin2006finite}. Such methods are well-suited for one-dimensional curvilinear boundary problems. Similarly, Costa et al. \cite{costa2022fft, wang2023fft} adopted a comparable approach, leveraging FFT to transform three-dimensional problems into two-dimensional elliptic problems for solving high-resolution turbulence. However, the applicability of these solvers is restricted to problems with one-dimensional curvilinear boundaries, limiting their generalization to more complex geometries.


Another primary direction of investigation involves non-body-conforming grid methods, exemplified by the Immersed Boundary Method (IBM)\cite{peskin1972flow}. Such methods construct interpolation schemes at the boundaries independently, thereby preserving the structure of the underlying regular grid. This characteristic makes them highly suitable for integration with FFT-based solution algorithms, leading to extensive research and numerous applications\cite{gilmanov2005hybrid, albuquerque2024numerical, chen2024synergistic, chen2024fourier, uhlmann2005immersed}. Subsequent work has focused on improving both efficiency and accuracy. For instance, Lian et al. \cite{lian2025overlapping} introduced an overlapping computation strategy to enhance parallel efficiency, while Zhao et al. \cite{zhao2023kernel} proposed kernel-free integration methods to reduce computational cost. In terms of accuracy, Ren et al. proposed a complex geometry handling method that combines the Matched Interface Boundary (MIB) technique with projection methods \cite{ren2022fft}, and subsequently extended it to fourth-order precision\cite{ren2023fft}, various types of boundary conditions \cite{li2023fast} and multi-physics problems\cite{song2024new}.

However, handling complex external boundary problems is notably challenging for IBM-based methods. Moreover, the absence of body-fitted grids makes it difficult to apply local mesh refinement, which is often crucial for the resolution of boundary layers. These limitations collectively restrict the applicability of non-body-conforming grid methods \cite{verzicco2023immersed}.

In this paper, we propose a novel FFT-based domain decomposition algorithm that overcomes some of the crucial limitations of the aforementioned methods. To the best of our knowledge, this specific approach has not been previously explored in the published literature. The proposed method can be conceptually divided into two primary components: a domain decomposition strategy for the global equation and an FFT-based GMRES \cite{saad1986gmres} algorithm for solving the resulting sub-equations. This algorithm relies on discrete equations expressed in matrix form. To ensure compatibility with this matrix formulation, an FFT-based solution algorithm, described through matrix eigenvalue transformations, is presented in this paper. In terms of computational performance, this algorithm exhibits the same $O\left(N \log N\right)$ complexity and inherent parallelism as solution algorithms based on spectral space expansion. Moreover, it can be more readily generalized to problems involving various types of boundary conditions.

In the domain decomposition algorithm proposed in this paper, the computational domain $\Omega$ is decomposed into multiple sub-regions $\Omega_i\left(i\in\mathbb{N}_+\right)$ with simple geometries, following the principles of DDMs. On each sub-region $\Omega_i$, sub-equations suitable for FFT-based algorithms can be directly constructed. The connections between sub-regions are described by the Schur complement.\cite{schur1917potenzreihen} Similar to general DDMs, the global equation is thus transformed into sub-systems of equations on each sub-region. 

The Schur complement disrupts the specific matrix structure required by FFT-based solvers, preventing a direct solution of the sub-equations. Therefore, this paper proposes an FFT-based preconditioned GMRES algorithm for solving these sub-equations. Specifically, the inverse of an approximate matrix that omits the Schur complement is employed as the preconditioner. This approach accelerates the Arnoldi process within GMRES by directly solving equations via FFT. While reducing the computational complexity of each iteration to $O\left(N \log N\right)$, this method also significantly decreases the number of iterations. The overall computational complexity is then evaluated through numerical results.

This paper is organized as follows. In Section \ref{sec:method}, we begin by reviewing the FFT algorithm for simple geometries based on matrix eigenvalue transformations, and then introduce the proposed FFT-based domain decomposition algorithm, including the FFT-preconditioned GMRES algorithm. To validate the solution efficiency of our proposed algorithm, Section \ref{sec:results} elaborates on the numerical results for the 2D Poisson equation within a cross-shape domain, including an analysis of the impact of parameters. Finally, Section \ref{sec:conclusion} provides a summary of the methods proposed in this study.

\section{Methodology}
\label{sec:method}
\subsection{FFT-based Poisson algorithm}
\label{ssec:FFT}
In this section, we begin by introducing the FFT-based Poisson equation solver algorithm for simple geometries, which forms the foundation of our proposed algorithm. Unlike traditional FFT-based solvers that rely on spectral expansion, the FFT-based solver for simple geometries proposed in this paper primarily focuses on low-complexity feature transformations of matrix equations, including the implementation of distinct FFT-based Poisson equation solvers for various boundary conditions.

The two-dimensional Poisson equation $\nabla^2u=f$ in $\Omega\subset\mathbb{R}^2$ is discretized by a central difference scheme, such that
\begin{equation}
    \label{eq:Poisson_central}
    \dfrac{p_{i+1,j}-2p_{i,j}+p_{i-1,j}}{\Delta_x^2}+\dfrac{p_{i,j+1}-2p_{i,j}+p_{i,j+1}}{\Delta_y^2}=f_{i,j},
\end{equation}
where $\Delta_x$ and $\Delta_y$ are the spatial steps on $x$ and $y$ direction, respectively.

Traditional FFT-based Poisson equation solvers are typically applicable only when $\Omega$ is a simple geometry, such as a rectangle in a Cartesian coordinate system or an annulus in a polar coordinate system. To simplify the presentation, we describe the discretization of the Poisson equation on a two-dimensional rectangular domain $\Omega$ in a Cartesian coordinate system. A uniform mesh is employed, with the number of nodes in the $x$ and $y$ directions being $\left(m+2\right)$ and $\left(n+2\right)$ respectively, including the boundary nodes. The central difference scheme \eqref{eq:Poisson_central} for the interior points can thus be written as the $mn \times mn$ block tridiagonal matrix equation
\begin{equation}
    \label{eq:mat_Poisson}
    \left[\begin{matrix}
        \boldsymbol{A}_y & \delta_x\boldsymbol{I} \\
        \delta_x\boldsymbol{I} & \boldsymbol{A}_y & \ddots\\
        & \ddots & \ddots & \delta_x\boldsymbol{I}\\
        & & \delta_x\boldsymbol{I} & \boldsymbol{A}_y
    \end{matrix}\right]
    \left[\begin{matrix}
        \boldsymbol{p}_{1}\\
        \boldsymbol{p}_{2}\\
        \vdots\\
        \boldsymbol{p}_{m}\\
    \end{matrix}\right]
    =
    \left[\begin{matrix}
        \boldsymbol{f}_{1}\\
        \boldsymbol{f}_{2}\\
        \vdots\\
        \boldsymbol{f}_{m}\\
    \end{matrix}\right].
\end{equation}
In \eqref{eq:mat_Poisson}, $\delta_x=\Delta_x^{-2}$ and $\delta_y=\Delta_y^{-2}$. The vector $\boldsymbol{p}=\left[\begin{matrix}p_{i,1}&p_{i,2}&\cdots&p_{i,n}\end{matrix}\right]^{\rm T}$ contains the unknown values along the $i$-th grid line, and the vector$\boldsymbol{f}$ is defined analogously. Assuming that the Poisson equation has Dirichlet boundary conditions $\left.p\right|_{\partial\Omega}=0$ on all four boundaries of the rectangular domain, the $n\times n$ submatrix $\boldsymbol{A}_y$ is
\begin{equation}
    \label{eq:submat_Dirichlet}
    \boldsymbol{A}_y=\boldsymbol{A}_y^{\rm D}=
    \left[\begin{matrix}
        -2\left(\delta_x+\delta_y\right) & \delta_y \\
        \delta_y & -2\left(\delta_x+\delta_y\right) & \ddots\\
        & \ddots & \ddots & \delta_y\\
        & & \delta_y & -2\left(\delta_x+\delta_y\right)
    \end{matrix}\right].
\end{equation}



The matrix $\boldsymbol{A}_y^{\rm D}$ is a symmetric Toeplitz matrix whose eigenvectors form the basis of the Discrete Sine Transform (DST). This property is exploited by traditional FFT-based solvers. The solution is expanded in the basis of these eigenvectors along the $y$-direction:
\begin{equation}
    \label{eq:Poisson_spectral}
    p_{i,j}=\sum_{k=1}^{n}\hat{p}_{i,k}\sin\left(\frac{jk\pi}{n+1}\right),
\end{equation}
where $\hat{p}_{i,k}$ are the spectral coefficients. Substituting this expansion into the finite difference scheme decouples the system of equations. For each spectral mode $k$, we obtain an independent tridiagonal system for the coefficients $\hat{p}_{i,k}$ along the $x$-direction:
\begin{equation}
    \label{eq:Poisson_spectral_equ}
    \delta_x\hat{p}_{i-1,k} + \lambda_k \hat{p}_{i,k} + \delta_x\hat{p}_{i+1,k} = \hat{f}_{i,k},
\end{equation}
where $\hat{f}_{i,k}$ are the transformed source terms and the $k$-th eigenvalue $\lambda_k$ of the submatrix $\boldsymbol{A}_y^{\rm D}$ is defined by
\begin{equation}
    \label{eq:Poisson_spectral_eigen}
    \lambda_k = -2(\delta_x+\delta_y) + 2\delta_y\cos(\frac{k\pi}{n+1}).
\end{equation}
It can be observed that these equations are independent, allowing the global equation to be transformed into multiple independent sub-equations, which is particularly suitable for parallel computation. Furthermore, each transformed sub-equation \eqref{eq:Poisson_spectral_equ} possesses a tridiagonal structure, which can be efficiently solved using the Thomas algorithm with a complexity of $O(m)$. Moreover, \eqref{eq:Poisson_spectral} can be directly computed with $O(n \log n)$ complexity using the discrete sine transform accelerated by FFT. Therefore, FFT-based solvers have become one of the fastest algorithms for solving Poisson equations in simple geometries.

However, traditional FFT-based Poisson equation solvers relying on spectral expansion face challenges posed by potential aliasing issues. Additionally, for various boundary conditions, the process of constructing basis functions can be complex, or it may be difficult to devise suitable spectral expansion methods. Therefore, the primary idea behind the FFT-based Poisson equation solver algorithm adopted in this study is the eigenvalue decomposition of \eqref{eq:mat_Poisson}. The FFT algorithm is then utilized to accelerate this transformation process. This algorithm shares similar properties with spectral expansion methods, specifically its ability to transform a global equation into independent sub-equations, thereby making it suitable for parallel computing.

This section details the FFT solver based on eigendecomposition, using Dirichlet boundary conditions $\left.p\right|_{\partial\Omega}=0$ as an illustrative example. The eigenvalues $\lambda_j^{\rm D}$ of matrix $\boldsymbol{A}_y^{\rm D}$ are same as \eqref{eq:Poisson_spectral_eigen}. Let the $j$-th eigenvector be $\boldsymbol{q}^{\rm D}_{j}$, whose elements are
\begin{equation}
    \label{eq:eigvector_Diri}
    {q}^{\rm D}_{j,k}=\sqrt{\frac{2}{n+1}}\sin{\dfrac{jk\pi}{n+1}},\quad \text{for } k={1, \dots, n}.
\end{equation}
Define $\boldsymbol{Q}_{\rm D}$ as the matrix whose columns are the eigenvectors $\boldsymbol{q}^{\rm D}_{j}$. Consequently, it can be further shown that $\boldsymbol{\Lambda}_y^{\rm D}=\boldsymbol{Q}_{\rm D}\boldsymbol{A}_y^{\rm D}\boldsymbol{Q}_{\rm D}^{\rm T}$, where $\boldsymbol{\Lambda}_y$ is a diagonal matrix with $\lambda_j^{\rm D}$ as its diagonal elements. By applying an eigenvalue transformation, equation \eqref{eq:mat_Poisson} can be transformed into
\begin{equation}
    \label{eq:eigen_trans}
    \left[\begin{matrix}
        \boldsymbol{\Lambda}_y & \delta_x\boldsymbol{I} \\
        \delta_x\boldsymbol{I} & \boldsymbol{\Lambda}_y & \ddots\\
        & \ddots & \ddots & \delta_x\boldsymbol{I}\\
        & & \delta_x\boldsymbol{I} & \boldsymbol{\Lambda}_y
    \end{matrix}\right]
    \left[\begin{matrix}
        \hat{\boldsymbol{p}}_{1}\\
        \hat{\boldsymbol{p}}_{2}\\
        \vdots\\
        \hat{\boldsymbol{p}}_{n}\\
    \end{matrix}\right]
    =
    \left[\begin{matrix}
        \hat{\boldsymbol{f}}_{1}\\
        \hat{\boldsymbol{f}}_{2}\\
        \vdots\\
        \hat{\boldsymbol{f}}_{n}\\
    \end{matrix}\right],
\end{equation}
where the vector $\hat{\boldsymbol{p}}_i=\boldsymbol{Q}^{\rm T}\boldsymbol{p}_i$ and $\hat{\boldsymbol{f}}_i=\boldsymbol{Q}^{\rm T}\boldsymbol{f}_i$.\par
$\hat{\boldsymbol{p}}'_j$ is defined by reordering elements of $\hat{\boldsymbol{p}}_j$, such that
\begin{equation}
    \label{eq:trans_order}
    \hat{\boldsymbol{p}}'_j=\left[\begin{matrix}
        \hat{u}_{1,j} & \hat{u}_{2,j} & \cdots & \hat{u}_{n,j}
    \end{matrix}\right]^{\rm T},
\end{equation}
and the vector $\hat{\boldsymbol{f}}'_j$ is defined analogously. \eqref{eq:eigen_trans} can be further rearranged into
\begin{equation}
    \label{eq:trans_eq}
    \left[\begin{matrix}
        \hat{\boldsymbol{A}}_1\\
        &\hat{\boldsymbol{A}}_2\\
        &&\ddots\\
        &&&\hat{\boldsymbol{A}}_m
    \end{matrix}\right]
    \left[\begin{matrix}
        \hat{\boldsymbol{p}}'_1\\
        \hat{\boldsymbol{p}}'_2\\
        \vdots\\
        \hat{\boldsymbol{p}}'_m\\
    \end{matrix}\right]
    =\left[\begin{matrix}
        \hat{\boldsymbol{f}}'_1\\
        \hat{\boldsymbol{f}}'_2\\
        \vdots\\
        \hat{\boldsymbol{f}}'_m\\
    \end{matrix}\right],
\end{equation}
where $\hat{\boldsymbol{A}}_k$ is the tridiagonal matrix forms that
\begin{equation}
    \label{eq:trans_subeq}
    \hat{\boldsymbol{A}}_k=
    \left[\begin{matrix}
        \lambda_k^{\rm D} & \delta_x \\
        \delta_x & \lambda_k^{\rm D} & \ddots \\
        & \ddots & \ddots & \delta_x \\
        & & \delta_x & \lambda_k^{\rm D}
    \end{matrix}\right].
\end{equation}
\eqref{eq:trans_eq} can be directly transformed into $n$ independent sub-equations, $\hat{\boldsymbol{A}}_k\hat{\boldsymbol{p}}'_k=\hat{\boldsymbol{f}}'_k$. Each of these sub-equations possesses a tridiagonal matrix on its left-hand side, as shown in \eqref{eq:trans_subeq}. This indicates that the FFT-based Poisson equation solver algorithm, when derived through eigenvalue transformation, shares the same characteristics as those obtained via spectral expansion.

The foundation of the FFT-based Poisson equation solver algorithm employed in this paper lies in the fast computation of $\boldsymbol{Q}_{\rm D}\boldsymbol{p}$ and $\boldsymbol{Q}_{\rm D}^{\rm T}\boldsymbol{p}$ for an arbitrary vector $\boldsymbol{p}$ of length $n$, achieved with $O(n \log n)$ complexity using the FFT algorithm. The core of this algorithm is to achieve the rapid computation of $\boldsymbol{Q}_{\rm D}\boldsymbol{p}$ and $\boldsymbol{Q}_{\rm D}^{\rm T}\boldsymbol{p}$ through a combination of the FFT algorithm and other low-complexity operations. It is straightforward to demonstrate that the DFT of a vector $\boldsymbol{p}$ is equivalent to the linear transformation
\begin{equation}
    \label{eq:Fourier_trans}
    \mathcal{F}\left[\boldsymbol{p}\right]=\boldsymbol{D}_n\boldsymbol{p},
\end{equation}
where $\boldsymbol{D}_n$ is a symmetric matrix of size $n$, and its $(j,k)$ element is given by
\begin{equation}
    \label{eq:Fourier_mat_element}
    d_{j,k}=e^{\frac{-2\pi{\rm i}\left(j-1\right)\left(k-1\right)}{n}},
\end{equation}
where ${\rm i}$ is the imaginary unit. The matrix $\boldsymbol{D}_n$ can be further decomposed into $\boldsymbol{D}_n=\boldsymbol{C}_n+{\rm i}\boldsymbol{S}_n$, where $\boldsymbol{C}_n$ and $\boldsymbol{S}_n$ represent its real and imaginary parts, respectively. It can be further shown that
\begin{equation}
    \label{eq:Dirichlet_mat_decomp}
    \boldsymbol{Q}_{\rm D}=\boldsymbol{K}_{\rm D}\boldsymbol{S}_{2n},
\end{equation}
where $\boldsymbol{K}_{\rm D}$ is an $n\times2n$ matrix with elements of 1 at each $(j, 2j)$ and zeros elsewhere. Therefore $\boldsymbol{Q}_{\rm D}\boldsymbol{p}_i$ can be computed with a complexity of $O(n \log n)$ using an FFT algorithm.\par

Although this section only presents the FFT-based method for 2D problems, the application of this method can be expanded to 3D cases. For 3D problems, another transformation is necessary to transform the 3D global equation into the independent 2D subequations.

\subsection{Boundary conditions}
\label{ssec:Bound}
The FFT-based solver algorithm presented previously, when applied to problems with Neumann and periodic boundary conditions, differs from the Dirichlet case only in the specific form of the transformation matrix $\boldsymbol{Q}$ and the corresponding fast transform algorithm.. This section introduces the low-complexity transform algorithms for these two types of boundaries, respectively.

For Neumann boundary conditions, where $\frac{\partial p}{\partial y}=0$ on the boundaries in the $y$-direction, the submatrix $\boldsymbol{A}_y$ forms
\begin{equation}
    \label{eq:submat_Neumann}
    \boldsymbol{A}_y=\boldsymbol{A}_y^{\rm N}=
    \left[\begin{matrix}
        -2\delta_x-\delta_y & \delta_y \\
        \delta_y & -2\left(\delta_x+\delta_y\right) & \ddots\\
        & \ddots & \ddots & \delta_y\\
        & & \delta_y & -2\delta_x-\delta_y
    \end{matrix}\right].
\end{equation}
The eigenvalues ${\lambda}^{\rm N}_{j}$ can be calculated as
\begin{equation}
    \label{eq:eigenvalues_Neumann}
    \lambda_j^{\rm N}=-2\left(\delta_x+\delta_y\right)+2\delta_y\cos\frac{(j-1)\pi}{n-1}, \quad \text{for } j=1, \dots, n.
\end{equation}
The matrix $\boldsymbol{Q}_{\rm N}$, containing the eigenvectors of $\boldsymbol{A}_y^{\rm N}$, is given by
\begin{equation}
    \label{eq:eigvector_Neumann_mat}
    \boldsymbol{Q}_{\rm N} = \bar{\boldsymbol{Q}}_{\rm N} + \frac{1+\sqrt{2}}{\sqrt{n}}\boldsymbol{E}_n,
\end{equation}
where $\boldsymbol{E}_n$ is the matrix where only the $n$-th column contains ones, with all other columns being zeros. The elements of matrix $\bar{\boldsymbol{Q}}_{\rm N}$ are given by
\begin{equation}
    \label{eq:eigvector_Neumann_bar}
    \bar{q}^{\rm N}_{j,k}=\sqrt{\frac{2}{n}}\sin{\dfrac{\left(jk-k-2nj\right)\pi}{2n}}
\end{equation}
Therefore, $\bar{\boldsymbol{Q}}_{\rm N}$ can be further decomposed as
\begin{equation}
    \label{eq:eigvector_Neumann_bar_decomp}
    \bar{\boldsymbol{Q}}_{\rm N}=\sqrt{\frac{2}{n}}\boldsymbol{K}_{\rm N}\left(\boldsymbol{S}_{n}\boldsymbol{\Delta}_{\rm C}+\boldsymbol{C}_{n}\boldsymbol{\Delta}_{\rm S}\right),
\end{equation}
where matrices $\boldsymbol{K}_{\rm N}$, $\boldsymbol{\Delta}_{\rm C}$, and $\boldsymbol{\Delta}_{\rm S}$ are all $n \times n$ diagonal matrices, with their $k$-th elements being $\left(-1\right)^k$, $\cos\frac{k\pi}{2n}$, and $\sin\frac{k\pi}{2n}$, respectively. Similar to the Dirichlet boundary conditions, when computing $\boldsymbol{Q}_{\rm N}\boldsymbol{p}$, the computational process can be progressively decomposed according to \eqref{eq:eigvector_Neumann_mat} and \eqref{eq:eigvector_Neumann_bar_decomp}. Since the left multiplication of vectors by matrices $\boldsymbol{C}_{n}$ and $\boldsymbol{S}_{n}$ can be calculated via FFT, and all other involved matrices are diagonal, the computational process, similar to the Dirichlet boundary case, can be transformed into one with a complexity of $O\left(n\log n\right)$.

It is worth noting that $\bar{\boldsymbol{Q}}_{\rm N}$ is not symmetric. Therefore, unlike the Dirichlet boundary conditions, when computing $\boldsymbol{Q}_{\rm N}^{\rm T}\boldsymbol{p}$, it is necessary to rely on the decomposition of $\boldsymbol{Q}_{\rm N}^{\rm T}$, i.e.,
\begin{equation}
    \label{eq:eigvector_Neumann_bar_transpose_decomp}
    \boldsymbol{Q}_{\rm N}^{\rm T}=\frac{1+\sqrt{2}}{\sqrt{n}}\boldsymbol{E}_n^{\rm T}+\sqrt{\frac{2}{n}}\left(\boldsymbol{\Delta}_{\rm C}\boldsymbol{S}_{n}+\boldsymbol{\Delta}_{\rm S}\boldsymbol{C}_{n}\right)\boldsymbol{K}_{\rm N}.
\end{equation}
This yields a computational process different from that for $\boldsymbol{Q}_{\rm N}\boldsymbol{p}$. However, it can be readily shown that this process also has a computational complexity of $O\left(n\log n\right)$.\par
In traditional spectral expansion methods, periodic boundary conditions typically the simplest to implement. However, for the algorithm proposed in this section, the matrix corresponding to periodic boundary conditions is more complex than the matrix decompositions for Dirichlet problems \eqref{eq:Dirichlet_mat_decomp} and Neumann boundary problems \eqref{eq:eigvector_Neumann_bar_decomp}. The matrix $\boldsymbol{A}_y$ corresponding to periodic boundary conditions is given by
\begin{equation}
    \label{eq:submat_Periodic}
    \boldsymbol{A}_y=\boldsymbol{A}_y^{\rm P}=
    \left[\begin{matrix}
        -2\left(\delta_x+\delta_y\right) & \delta_y &  & \delta_y\\
        \delta_y & -2\left(\delta_x+\delta_y\right) & \ddots\\
        & \ddots & \ddots & \delta_y\\
        \delta_y & & \delta_y & -2\left(\delta_x+\delta_y\right)
    \end{matrix}\right].
\end{equation}
The eigenvalues ${\lambda}^{\rm P}_{j}$ can be calculated as
\begin{equation}
    \label{eq:eigenvalues_Periodic_corrected}
    \lambda_j^{\rm P}=-2\left(\delta_x+\delta_y\right)+2\delta_y\cos\frac{2\pi(j-1)}{n}, \quad \text{for } j=1, \dots, n.
\end{equation}
The elements of the matrix $\boldsymbol{Q}_{\rm P}$, which is composed of the eigenvectors of $\boldsymbol{A}_y^{\rm P}$, are given by
\begin{equation}
    \label{eq:Q_periodic_elements}
    q_{j,k}^{\rm P}=\sqrt{\frac{2}{n}}\sin\left[\frac{2j\pi}{n}\left\lfloor\frac{n+1-j}{2}\right\rfloor+\frac{\pi\left(2k-1\right)}{4}\right]
\end{equation}
where $\left\lfloor\cdot\right\rfloor$ denotes the floor function. Therefore, $\boldsymbol{Q}_{\rm P}$ can be decomposed as
\begin{equation}
    \label{eq:eigvector_Periodic_decomp}
    \boldsymbol{Q}_{\rm P}=\sqrt{\frac{2}{n}}\left(\boldsymbol{S}^{\rm P}\boldsymbol{\Delta}_{\rm PC}+\boldsymbol{C}^{\rm P}\boldsymbol{\Delta}_{\rm PS}\right),
\end{equation}
where matrices $\boldsymbol{\Delta}_{\rm PC}$ and $\boldsymbol{\Delta}_{\rm PS}$ are diagonal matrices whose $k$-th diagonal elements are $\cos\frac{\left(2k-1\right)\pi}{4}$ and $\sin\frac{\left(2k-1\right)\pi}{4}$, respectively. The decomposition in Eq. \eqref{eq:eigvector_Periodic_decomp} holds if and only if the product $n$ is an even number. To satisfy this condition, we assume throughout this work that $n$ is an even integer. Accordingly, the elements of matrix $\boldsymbol{S}^{\rm P}$ are
\begin{equation}
    \label{eq:Periodic_mat_s_p_element}
    s^{\rm P}_{j,k}=\left\{\begin{aligned}
        &\sin\frac{j\pi\left(n-2l\right)}{n}, &l\in\mathbb{N},k=2l\\
        &\sin\frac{j\pi\left(n-2l+2\right)}{n}. &l\in\mathbb{N},k=2l+1
    \end{aligned}\right.
\end{equation}
Similarly, the elements of $\boldsymbol{C}^{\rm P}$ are
\begin{equation}
    \label{eq:Periodic_mat_c_p_element}
    c^{\rm P}_{j,k}=\left\{\begin{aligned}
        &\cos\frac{j\pi\left(n-2l\right)}{n}, &l\in\mathbb{N},k=2l\\
        &\cos\frac{j\pi\left(n-2l+2\right)}{n}. &l\in\mathbb{N},k=2l+1
    \end{aligned}\right.
\end{equation}
To facilitate an FFT-based implementation, we partition the matrices $\boldsymbol{C}^{\rm P}$ and $\boldsymbol{S}^{\rm P}$ based on column parity into four intermediate $n/2 \times n$ matrices. The matrix $C^{\rm PE}$ is constructed from the even-indexed columns of $\boldsymbol{C}^{\rm P}$, and $C^{\rm PO}$ from the odd-indexed columns. Their elements are given by
\begin{equation}
    \label{eq:Periodic_mat_c_pe_element}
    c^{\rm PE}_{j,l}=\cos\frac{j\pi\left(n-2l\right)}{n}=\left(-1\right)^j\cos\frac{2\pi jl}{n},
\end{equation}
and
\begin{equation}
    \label{eq:Periodic_mat_c_po_element}
    c^{\rm PO}_{j,l}=\cos\frac{j\pi\left(n-2l+2\right)}{n}.
\end{equation}

The matrices $S^{\rm PE}$ and $S^{\rm PO}$ are constructed analogously from the even and odd columns of $\boldsymbol{S}^{\rm P}$, respectively. It can be shown that $C^{\rm PE}$ and $S^{\rm PE}$ correspond to the upper halves of $\boldsymbol{K}_{\rm N}S_{2n}$ and $\boldsymbol{K}_{\rm N}S_{2n}$, whose low-complecity algorithm has been introduced above. Consequently, the matrix-vector product $\boldsymbol{Q}_{\rm P}\boldsymbol{p}$ can be efficiently computed using two FFTs, following a strategy analogous to that employed for the Neumann problem.

In this section, we have detailed an FFT-based algorithm derived from the eigendecomposition of the discrete operator. A notable feature of this method is that, contrary to classic FFT-based algorithm where periodic problems are simplest. Nevertheless, we have established a low-complexity, FFT-based solver for each of the three boundary condition types (Dirichlet, Neumann, and periodic). This unified framework is fundamental to the domain decomposition algorithm presented in the subsequent sections.


\subsection{Domain decomposition framework}
\label{ssec:ddm}
\begin{figure}[H]
    \centering
    \includegraphics[width=0.4\textwidth]{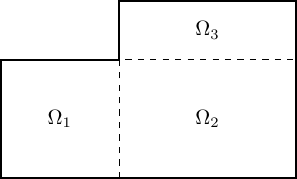}
    \caption{The schematic of the L-shaped region $\Omega$ with subdomains $\Omega_1$, $\Omega_2$, and $\Omega_3$.}
    \label{fig:concat}
\end{figure}

This section introduces the proposed DDM, using the two-dimensional L-shaped region shown in Fig.\ref{fig:concat} as an example. The two-dimensional L-shaped region in Fig.\ref{fig:concat} can be decomposed into three connected rectangular sub-regions, labeled as $\Omega_1$, $\Omega_2$, and $\Omega_3$, respectively.

The global equation has the block matrix structure, which forms
\begin{equation}
    \label{eq:block_global}
    \left[\begin{matrix}
        \boldsymbol{A}_1 & \boldsymbol{R}_{12} & \\
        \boldsymbol{R}_{21} & \boldsymbol{A}_2 & \boldsymbol{R}_{23} \\
        & \boldsymbol{R}_{32} & \boldsymbol{A}_3 \\
    \end{matrix}\right]
    \left[\begin{matrix}
        \boldsymbol{p}_{1}\\
        \boldsymbol{p}_{2}\\
        \boldsymbol{p}_{3}
    \end{matrix}\right]
    =
    \left[\begin{matrix}
        \boldsymbol{f}_{1}\\
        \boldsymbol{f}_{2}\\
        \boldsymbol{f}_{3}
    \end{matrix}\right],
\end{equation}
where the sparse matrices $\boldsymbol{R}_{ij}$ represent the coupling across the interface between subdomains $\Omega_i$ and $\Omega_j$. For the self-adjoint Poisson operator, these matrices satisfy $\boldsymbol{R}_{ij}=\boldsymbol{R}_{ji}^{\rm T}$. Each sub-equation $\boldsymbol{A}_i\boldsymbol{q}_i=\boldsymbol{f}_i\left(i=1,2,3\right)$ can be directly interpreted as a Poisson equation within the sub-domain $\Omega_i$, where the external boundary conditions are consistent with those of the original domain $\Omega$, and the boundary conditions at the interfaces between sub-domains are treated as Dirichlet boundary conditions.

Therefore, each sub-equation possesses at least one Dirichlet boundary, which implies that $\boldsymbol{A}_1$, $\boldsymbol{A}_2$, and $\boldsymbol{A}_3$ are all invertible matrices. Consequently, by applying the Schur complement, the global system \eqref{eq:block_global} can be transformed into three sub-equations:
\begin{subequations}
    \label{eq:Schur_subeq}
    \begin{align}
        \boldsymbol{p}_{1}&=\boldsymbol{q}_{1}-\boldsymbol{A}_{1}^{-1}\boldsymbol{R}_{12}\boldsymbol{p}_{2},\label{eq:Schur_subeq_1}\\
        \left(\boldsymbol{A}_{2}-\boldsymbol{S}_{21}-\boldsymbol{S}_{23}\right)\boldsymbol{p}_{2}&=\boldsymbol{f}_{2}-\boldsymbol{R}_{21}\boldsymbol{q}_{1}-\boldsymbol{R}_{23}\boldsymbol{q}_{3},\label{eq:Schur_subeq_2}\\
        \boldsymbol{p}_{3}&=\boldsymbol{q}_{3}-\boldsymbol{A}_{3}^{-1}\boldsymbol{R}_{32}\boldsymbol{p}_{2},\label{eq:Schur_subeq_3}
    \end{align}
\end{subequations}
where $\boldsymbol{S}_{ij}=\boldsymbol{R}_{ij}\boldsymbol{A}_{i}^{-1}\boldsymbol{R}_{ji}$ and $\boldsymbol{q}_{i}=\boldsymbol{A}_{i}^{-1}\boldsymbol{f}_{i}$.
 
Since each subdomain is rectangular, the matrices $\boldsymbol{A}_1$, $\boldsymbol{A}_2$, and $\boldsymbol{A}_3$ all possess the block-tridiagonal structure shown in \eqref{eq:mat_Poisson}, which is required for the FFT-based solver. Consequently, the equations $\boldsymbol{q}_{i}=\boldsymbol{A}_{i}^{-1}\boldsymbol{f}_{i}$ can be solved efficiently using the algorithm described in Section \ref{ssec:FFT}. However, it is evident that the structure of the left-hand side matrix $\boldsymbol{L}_{2}=\boldsymbol{A}_{2}-\boldsymbol{S}_{21}-\boldsymbol{S}_{23}$ in \eqref{eq:Schur_subeq_2} does not satisfy the conditions for direct solution using FFT-based algorithms. Therefore, an iterative solution algorithm is proposed in this paper.

\subsection{FFT-based GMRES algorithm}
\label{ssec:GMRES}
Let equation \eqref{eq:Schur_subeq_2} be simplified as $\boldsymbol{L}_{2}\boldsymbol{p}_{2}=\boldsymbol{f}'_{2}$, where $\boldsymbol{f}'_{2} = \boldsymbol{f}_{2}-\boldsymbol{R}_{21}\boldsymbol{q}_{1}-\boldsymbol{R}_{23}\boldsymbol{q}_{3}$. Due to the sparsity of $\boldsymbol{S}_{21}$ and $\boldsymbol{S}_{23}$, a fixed-point iteration can be directly formulated as
\begin{equation}
    \label{eq:basic_iter}
    \boldsymbol{p}^{\left(k+1\right)}=\boldsymbol{A}_2^{-1}\left(\boldsymbol{S}_{21}+\boldsymbol{S}_{23}\right)\boldsymbol{p}^{\left(k\right)}+\boldsymbol{q}'_{2},
\end{equation}
where $\boldsymbol{p}^{\left(k\right)}$ represents the solution at the $k$-th iteration step. In this iterative process, each iteration step effectively involves solving the equation
\begin{equation}
    \label{eq:basic_iter_equ}
    \boldsymbol{A}_2\boldsymbol{p}^{\left(k+1\right)}=\boldsymbol{S}_{21}+\boldsymbol{S}_{23}\boldsymbol{p}^{\left(k\right)}+\boldsymbol{f}'_{2}.
\end{equation}
This equation can be directly solved using the algorithm described in Section \ref{ssec:FFT}, which theoretically yields an $O\left(n\log n\right)$ complexity for each iteration step. However, a critical issue is that the matrix $\boldsymbol{A}_2$ is not diagonally dominant. Consequently, the convergence of this simple fixed-point iteration is not guaranteed. To address this problem, an improved Generalized Minimal Residual (GMRES) algorithm, enhanced by FFT, is proposed in this paper to efficiently solve \eqref{eq:Schur_subeq_2} based on the matrix structure.

Since $\boldsymbol{L}_{2}=\left(\boldsymbol{A}_{2}-\boldsymbol{S}_{21}-\boldsymbol{S}_{23}\right)$ is generally not symmetric, algorithms like Conjugate Gradient (CG), which require specific matrix structures, are not applicable. GMRES is therefore one of the most widely employed choices, balancing both convergence efficiency and stability. The computational process of the GMRES algorithm can be decomposed into the Arnoldi process and a least-squares problem on the Krylov subspace. Let the dimension of the Krylov subspace be $m$. Assuming the initial residual of the iteration is $\boldsymbol{r}^{\left(0\right)}$, the iterative process of the GMRES algorithm is as follows:

\begin{enumerate}
    \item Generate the Krylov subspace basis vectors via $\boldsymbol{w}_{i}=\boldsymbol{L}_{2}\boldsymbol{w}_{i-1}\left(i=2,\cdots,m\right)$, where $\boldsymbol{w}_{1}=\boldsymbol{r}^{\left(0\right)}$, .
    \item Orthogonalize the Krylov subspace basis vectors $\boldsymbol{w}_{i}$ using the Gram-Schmidt method and achieve $V_d = [\boldsymbol{v}_1, \dots, \boldsymbol{v}_d]$.
    \item Solve the least-squares problem $y_d^* = \arg \min_{y_d \in \mathbb{R}^d} \|\left\|\boldsymbol{r}^{\left(0\right)}\right\|_2 \boldsymbol{e}_1 - \bar{H}_d \boldsymbol{y}_d\|_2$  on the constructed Krylov space, where $\bar{H}_d$ is the Hessenberg matrix obtained from the Arnoldi process.
    \item Update the approximate solution: Calculate $\boldsymbol{p}^{\left(k+1\right)} = \boldsymbol{p}^{\left(k\right)} + V_d \boldsymbol{y}_d^*$.
    \item If the residual $\boldsymbol{r}^{\left(k\right)}$ falls below a preset tolerance, the iteration terminates; otherwise, the process continues.
\end{enumerate}

It is straightforward to demonstrate that for general problems, the dimension of the Krylov subspace, $m$, is typically chosen to be much smaller than the number of unknowns in the subdomain, $n_2$. Consequently, the Arnoldi process often constitutes the most time-consuming part of the GMRES solver. Therefore, the core idea of the optimized GMRES algorithm proposed in this paper is to use $\boldsymbol{A}_{2}$ as a preconditioner is to employ $\boldsymbol{A}_{2}^{-1}$ as a right preconditioner. This transforms the system and accelerates convergence, while the FFT-based solver for $\boldsymbol{A}_2$ allows the preconditioning step to be performed efficiently.\par
Specifically, in the proposed method, the equation $\boldsymbol{L}_{2}\boldsymbol{p}_{2}=\boldsymbol{f}'_{2}$ is transformed into
\begin{equation}
    \label{eq:precond_iter}
    \boldsymbol{p}_{2}-\boldsymbol{A}_{2}^{-1}\left(\boldsymbol{S}_{21}+\boldsymbol{S}_{23}\right)\boldsymbol{p}_{2}=\boldsymbol{A}_{2}^{-1}\boldsymbol{f}'_{2}.
\end{equation}
The generation of Krylov subspace basis vectors, originally calculated by $\boldsymbol{w}_{i}=\boldsymbol{L}_{2}\boldsymbol{w}_{i-1}$, is replaced by
\begin{equation}
    \label{eq:precond_basis}
    \boldsymbol{w}_{i}=\boldsymbol{w}_{i-1}-\boldsymbol{w}_{i-1}'
\end{equation}
where $\boldsymbol{w}_{i-1}'$ is obtained by solving the equation
\begin{equation}
    \label{eq:precond_basis_equ}
    \boldsymbol{A}_{2}\boldsymbol{w}_{i-1}'=\left(\boldsymbol{S}_{21}+\boldsymbol{S}_{23}\right)\boldsymbol{w}_{i-1}
\end{equation}
using the FFT-based algorithm introduced in \ref{ssec:FFT}.

In summary, the proposed FFT-based domain decomposition algorithm for solving the global system in Eq. \eqref{eq:block_global} proceeds as \ref{alg:fft_dd_solver}.

\begin{algorithm}[H] 
\caption{FFT-based DDM algorithm}
\label{alg:fft_dd_solver}
    Solve $\boldsymbol{A}_{1}\boldsymbol{q}_{1}=\boldsymbol{f}_{1}$ for $\boldsymbol{q}_{1}$ using the FFT-based solver.\\
    Solve $\boldsymbol{A}_{3}\boldsymbol{q}_{3}=\boldsymbol{f}_{3}$ for $\boldsymbol{q}_{3}$ using the FFT-based solver.\\
    Form the modified right-hand side $\boldsymbol{f}'_{2} = \boldsymbol{f}_{2}-\boldsymbol{R}_{21}\boldsymbol{q}_{1}-\boldsymbol{R}_{23}\boldsymbol{q}_{3}$.\\
    Solve the preconditioned system from Eq. \eqref{eq:precond_iter} for $\boldsymbol{p}_{2}$ using the FFT-based GMRES algorithm.\\
    Form the right-hand sides $\boldsymbol{f}'_1 = \boldsymbol{f}_{1}-\boldsymbol{R}_{12}\boldsymbol{p}_{2}$ and $\boldsymbol{f}'_3 = \boldsymbol{f}_{3}-\boldsymbol{R}_{32}\boldsymbol{p}_{2}$.\\
    Solve $\boldsymbol{A}_{1}\boldsymbol{p}_{1}=\boldsymbol{f}'_1$ for $\boldsymbol{p}_{1}$ using the FFT-based solver.\\
    Solve $\boldsymbol{A}_{3}\boldsymbol{p}_{3}=\boldsymbol{f}'_3$ for $\boldsymbol{p}_{3}$ using the FFT-based solver.
\end{algorithm}

This constitutes the general procedure for the FFT-based domain decomposition algorithm proposed in this paper. This process can also be directly extended to solve Poisson equations within more complex concatenated geometries than that shown in Figure \ref{fig:concat}. The computational procedure of the entire domain decomposition algorithm categorizes the sub-domains into two types: As exemplified by sub-domains $\Omega_1$ and $\Omega_3$ in Figure \ref{fig:concat}, those sub-domains that are connected to only one other sub-domain can be directly solved using the FFT-based algorithm introduced in Section \ref{ssec:FFT}; these are thus termed independent sub-domains. Conversely, sub-domains like $\Omega_2$, which are adjacent to multiple other sub-domains, require their sub-equations to be solved using the proposed GMRES algorithm; these are referred to as coupled sub-domains. By employing the two distinct algorithms proposed in this paper, fast solutions with $O\left(n\log n\right)$ complexity are achieved for both types of sub-domains. It is worth noting that the subspace dimension $m$ plays a crucial role in determining the convergence rate of the iteration. Consequently, the overall computational complexity of the solution algorithm is jointly determined by the equation size and the subspace dimension $m$. This paper provides a detailed evaluation of the overall complexity of the solution algorithm based on numerical results in Section \ref{ssec:effi}.

\section{Numerical results}
\label{sec:results}
In this research, all programs used for algorithm evaluation are in-house codes, developed in C++ and utilizing the highly optimized FFTW library for all FFT computations \cite{FFTW05, JohnsonFrigo07:splitradix}. All source code, including the specific versions employed in this research, is open-source and publicly available on \url{https://github.com/ffskibkwi/Focalors_Poisson}.

Our solver provides two distinct parallel implementations: one based on OpenMP for shared-memory environments, and another based on MPI for distributed-memory systems. Both parallel versions are made available in the aforementioned open-source repository. For consistency and relevance to single-node performance, all computational performance results, including execution times and scalability metrics presented in this study, are exclusively derived from the OpenMP-parallelized version of the code.

All computations were executed on a single compute node, configured with two 26-core Intel Xeon CPUs (total 52 cores) operating at a base frequency of 2.1 GHz, and equipped with 256 GB of node memory.

\subsection{Precision validation}

\begin{figure}[H]
    \centering
    \includegraphics[width=5.5cm]{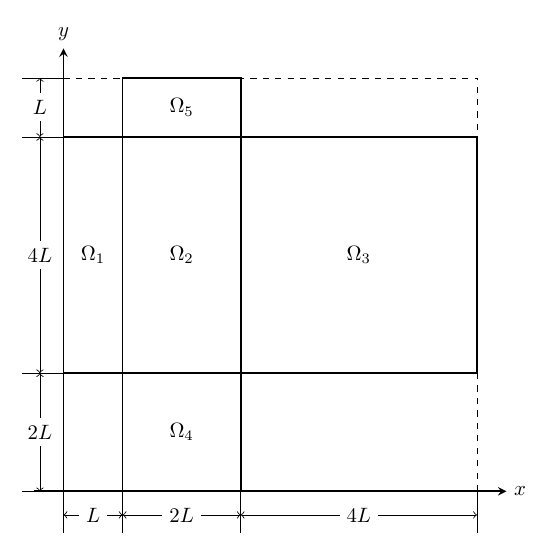}
    \caption{Cross-shape geometry.}
    \label{fig:cross_geometry}
\end{figure}

In this study, to evaluate the reliability and computational efficiency of the proposed Poisson equation solver, we construct a test case on a two-dimensional cross-shaped domain, as depicted in Figure \ref{fig:cross_geometry}.

In Figure \ref{fig:cross_geometry}, we defined the scale factor as $L=1/7$, and the equations were set with the following boundary conditions:
\begin{equation}
    \label{eq:cross_equation}
    \left\{\begin{aligned}
    &\dfrac{\partial^2p}{\partial x^2}+\dfrac{\partial^2p}{\partial y^2}=f\left(x,y\right),\\
    &\left.p\right|_{x=0}=\left.p\right|_{y=0}=\left.p\right|_{x=7L}=\left.p\right|_{y=7L}=0,\\
    &\left.\dfrac{\partial p}{\partial x}\right|_{x=L}=\left.\dfrac{\partial p}{\partial y}\right|_{y=2L}=\left.\dfrac{\partial p}{\partial x}\right|_{x=3L}=\left.\dfrac{\partial p}{\partial y}\right|_{y=6L}=0,\\
    \end{aligned}\right.
\end{equation}
This indicates that Dirichlet boundary conditions with a value of zero are applied on the outer edges, while Neumann boundary conditions with zero normal gradient are applied on the other edges. It can be demonstrated that the Poisson equation has the following analytical solution:

\begin{equation}
    \label{eq:cross_analytical_solultion}
    p\left(x,y\right)=\sin\left[\psi_x\left(x\right)\right]\sin\left[\psi_y\left(y\right)\right],
\end{equation}
where the functions $\psi_x(x)$ and $\psi_y(y)$ are constructed as piecewise polynomials to satisfy the boundary conditions, i.e.,
\begin{equation}
    \label{eq:cross_ana_sol_psi}
    \begin{aligned}
        &\psi_x\left(x\right)=x\left[A_1+A_2\left(x-L\right)+A_3\left(x^2-4Lx+3L^2\right)\right],\\
        &A_1=\dfrac{\pi}{2L},\quad A_2=\dfrac{1}{3L}\left(\dfrac{\pi}{2L}-A_1\right),\quad A_3=\dfrac{1}{7L}\left[-\dfrac{\pi}{8L^2}-A_2\right],\\
        &\psi_y\left(y\right)=y\left[B_1+B_2\left(y-2L\right)+B_3\left(y^2-8Ly+12L^2\right)\right],\\
        &B_1=\dfrac{\pi}{4L},\quad B_2=\dfrac{1}{6L}\left(\dfrac{\pi}{4L}-B_1\right),\quad B_3=\dfrac{1}{7L}\left[-\dfrac{\pi}{4L^2}-B_2\right].
    \end{aligned}
\end{equation}
The righthand term $f\left(x,y\right)$ is defineded by the Laplacian of the analytical solution $p\left(x,y\right)$.\par

On a uniform grid, the spatial step size is set as $\Delta x=\Delta y=L/k_n$, where $k_n$ is a given constant controlling the number of grid points. By comparing with the analytical solution, we evaluated the spatial accuracy order of the algorithm.

\begin{figure}[H]
    \centering
    \includegraphics[width=7.5cm]{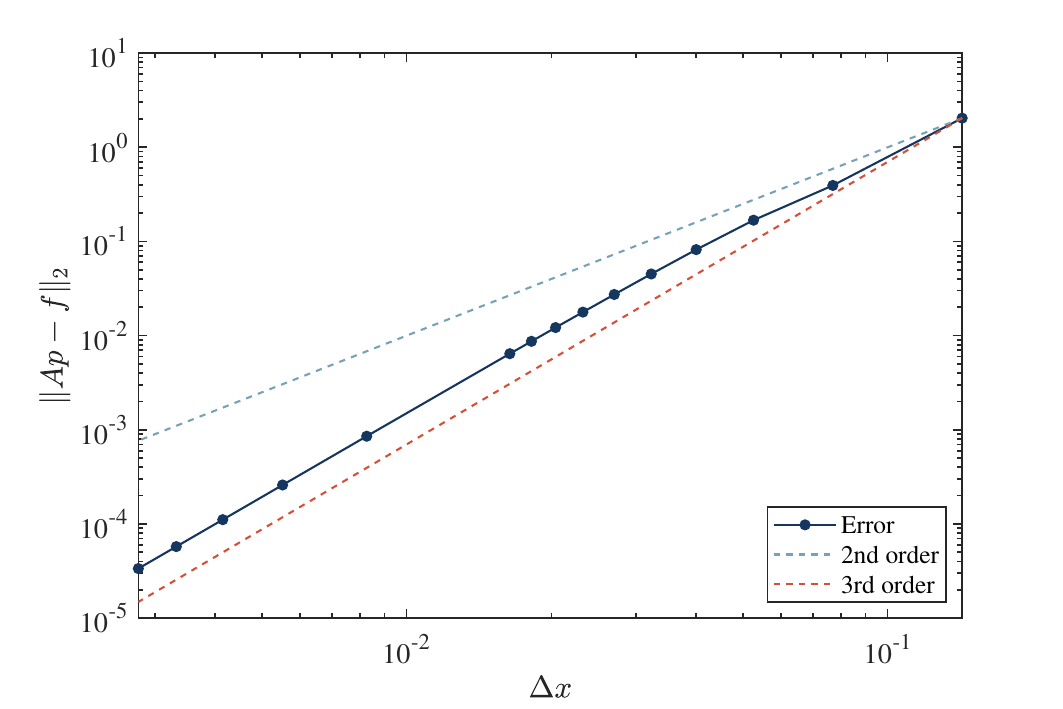}
    \caption{Precision order of cross-shape geometry.}
    \label{fig:cross_precision}
\end{figure}

As illustrated in Figure \ref{fig:cross_precision}, due to the second-order central difference scheme, the solver achieves the expected second-order rate of convergence, confirming that the proposed domain decomposition algorithm correctly reproduces the precision of the underlying difference scheme.

\subsection{GMRES iteration}

\begin{figure}[H]
	\centering
	\subfloat[\centering]{\includegraphics[width=5.5cm]{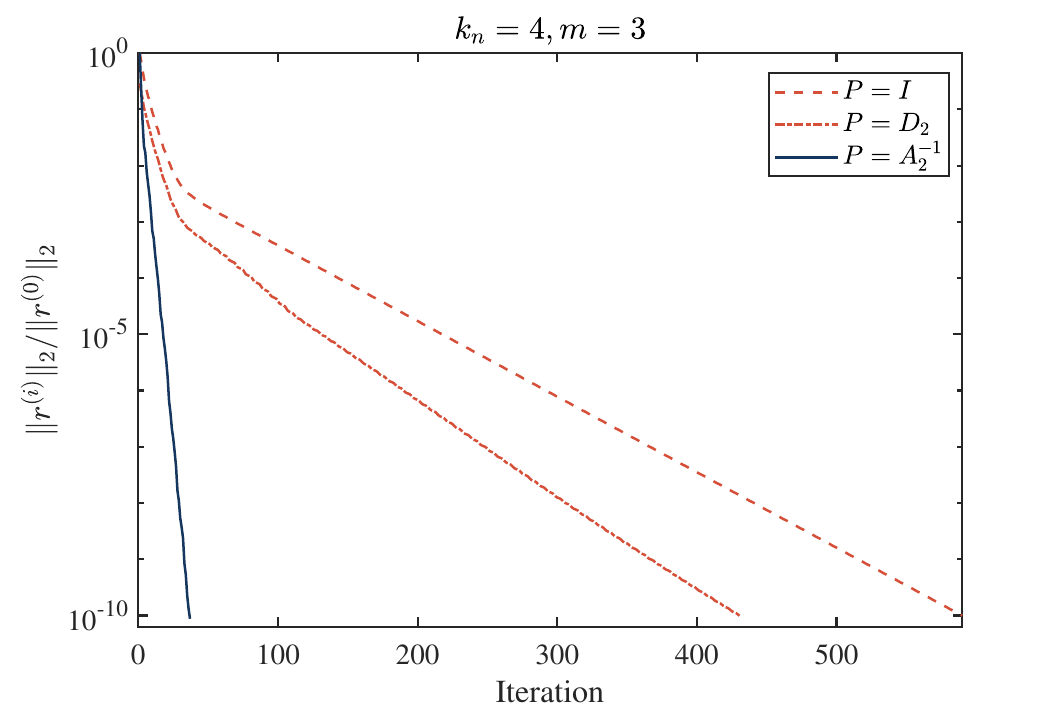}}
	\subfloat[\centering]{\includegraphics[width=5.5cm]{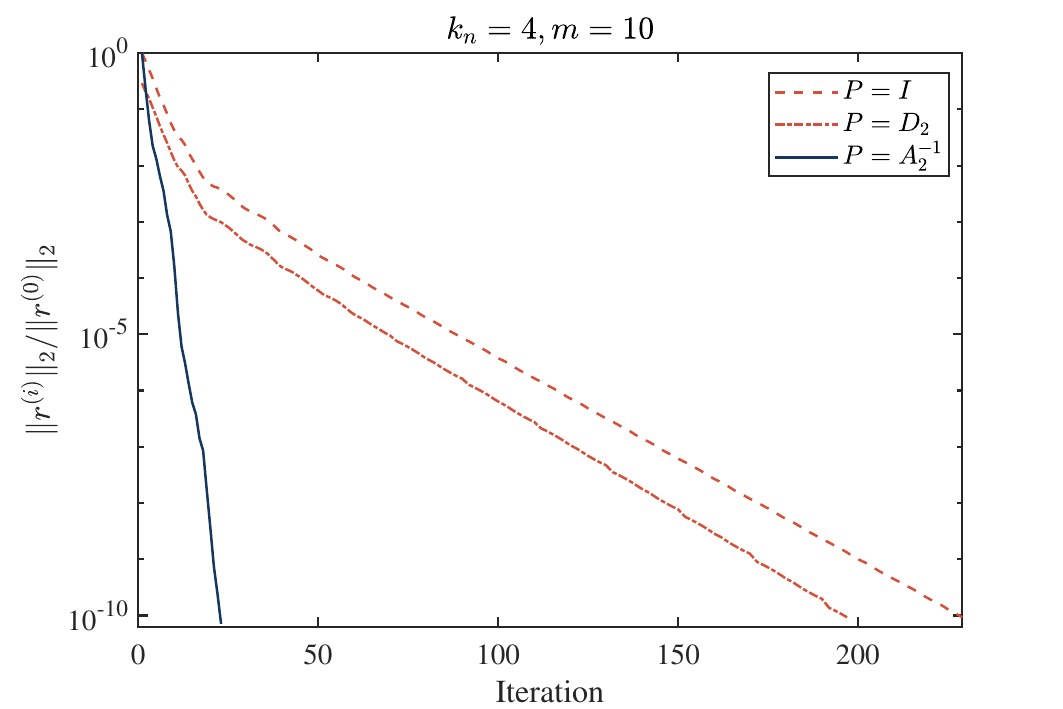}}

	\subfloat[\centering]{\includegraphics[width=5.5cm]{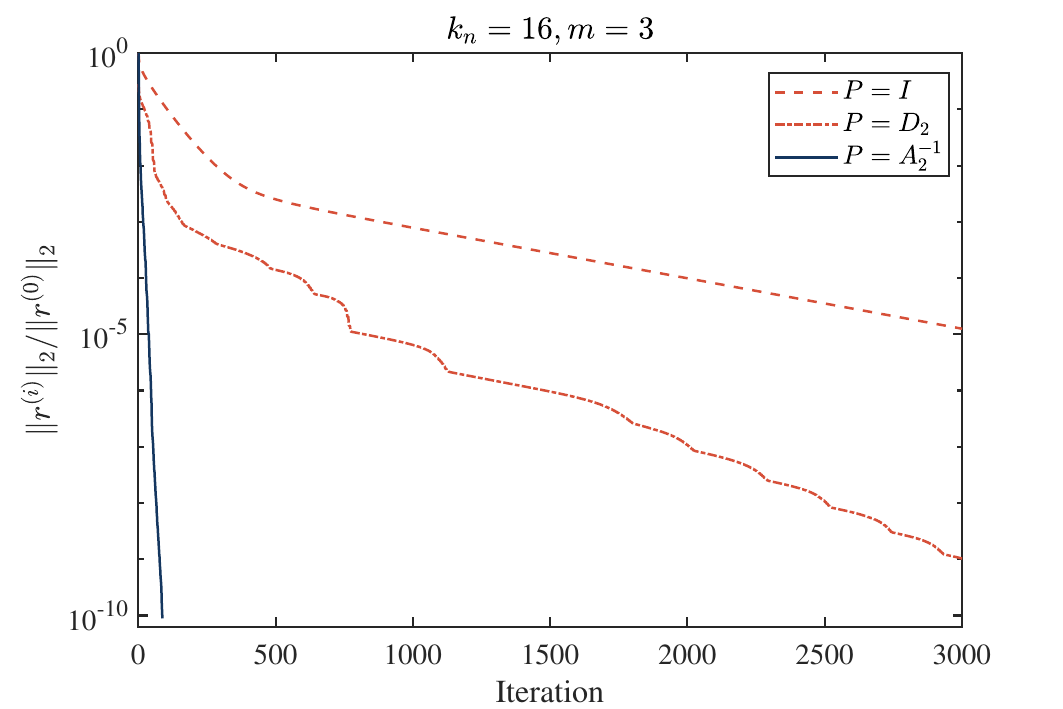}}
	\subfloat[\centering]{\includegraphics[width=5.5cm]{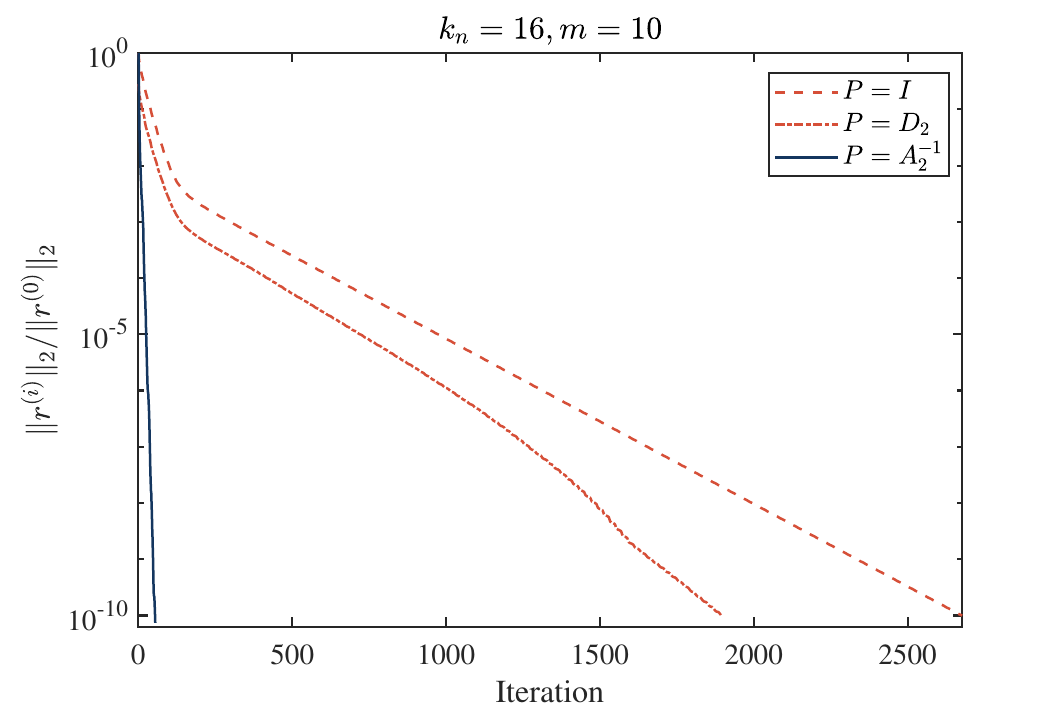}}
	\caption{Comparison among the proposed preconditioner and classic preconditioners with different dimensions of subspaces and equation sizes}
	\label{fig:precond_comp}
\end{figure}

As introduced in Section 2, the cross-flow field in Figure \ref{fig:cross_geometry} is decomposed into five rectangular subregions. As introduced in Section \ref{sec:method}, the cross-flow field in Figure \ref{fig:cross_geometry} is decomposed into five rectangular subregions. The solution of the central subdomain $\Omega_2$ is achieved by solving the Schur complement system via the propoesed preconditioned GMRES iteration. To demonstrate the benefits of the proposed FFT-based preconditioner in enhancing iterative performance, we solve the sub-equation on $\Omega_2$ using different preconditioners.  We consider two subspace dimensions $m$ and two equation sizes $k_n$ and evaluate the residual reduction during the iterative process. The results are presented in Figure \ref{fig:precond_comp}. To demonstrate the benefits of the proposed FFT-based preconditioner in enhancing iterative performance, we solve the sub-equation on $\Omega_2$ using different preconditioners.  We consider two subspace dimensions $m$ and two equation sizes $k_n$ and evaluate the residual reduction during the iterative process. The results are presented in Figure \ref{fig:precond_comp}.

In the comparison presented in Figure \ref{fig:precond_comp}, we evaluate the acceleration in convergence achieved by two preconditioning operators. The case $P=I$ indicates that no preconditioner is applied, and the sub-equation with the Schur complement, $(A_2-S_2)\boldsymbol{p}_2=\boldsymbol{f}_2'$, is solved directly. $P=D_2$ represents the use of the Jacobi preconditioner, a widely adopted classical preconditioner. Its principle involves using the inverse of the matrix $D_2$, formed by the main diagonal elements of $(A_2-S_2)$, as the preconditioning operator. $P=A_2^{-1}$ represents the proposed preconditioning operator as described in Section \ref{ssec:GMRES}. As can be observed, the proposed preconditioner not only guarantees convergence for all equation sizes tested but also exhibits a significantly higher convergence rate than the two methods used for comparison. The number of iterations required to reach a given tolerance is reduced by a factor of more than seven. This reduction in iteration count is fundamental to the efficiency of the proposed algorithm. By employing a low-complexity computation for each iteration, the proposed algorithm achieves a notable reduction in the total number of iterations.\par

\begin{figure}[H]
	\centering
	\subfloat[\centering]{\includegraphics[width=5.5cm]{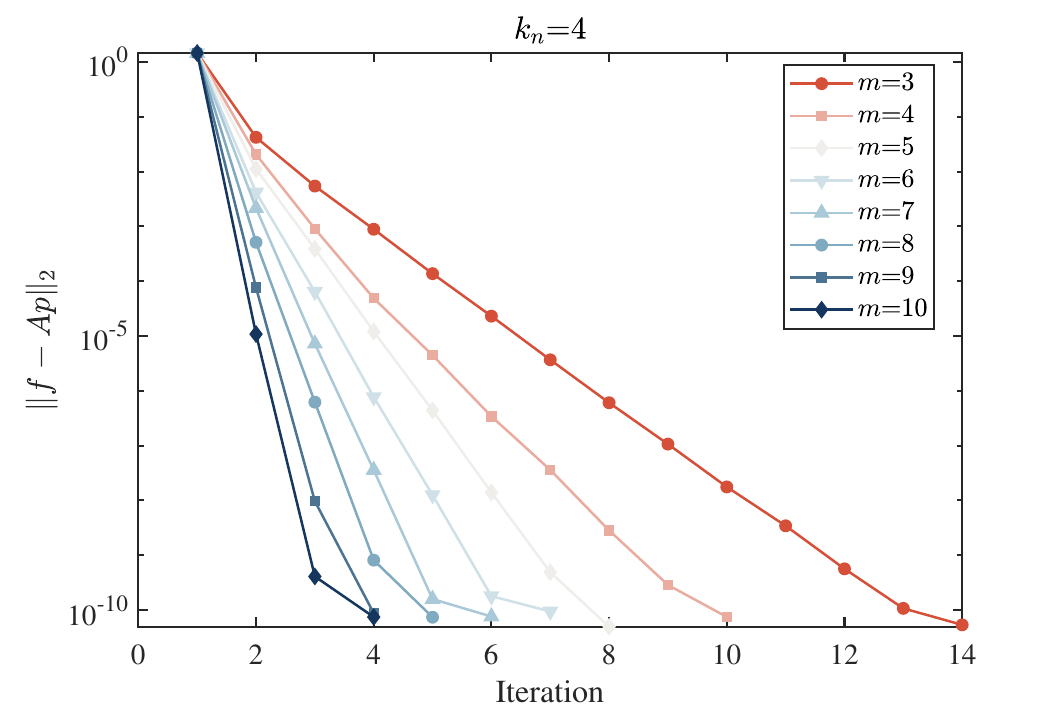}}
	\subfloat[\centering]{\includegraphics[width=5.5cm]{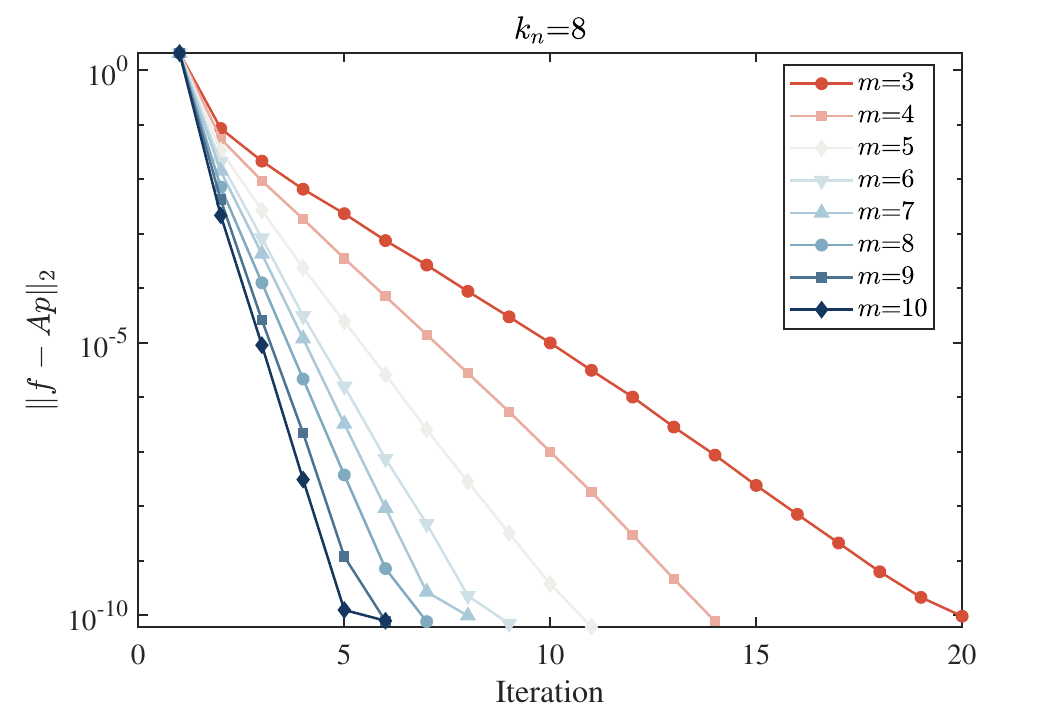}}

	\subfloat[\centering]{\includegraphics[width=5.5cm]{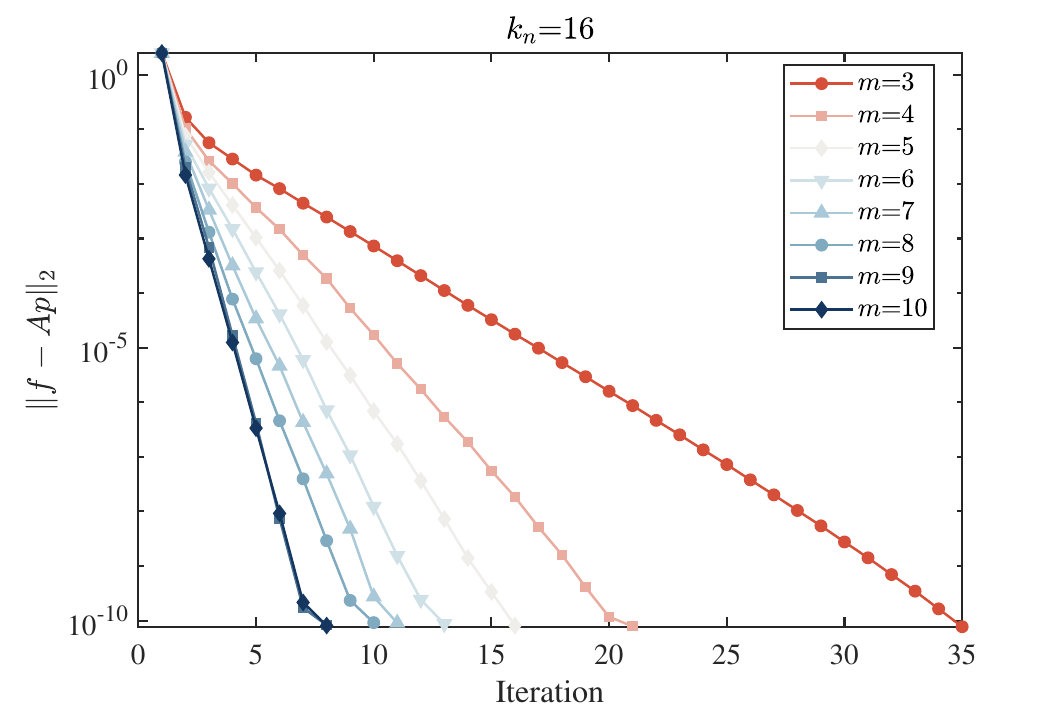}}
	\subfloat[\centering]{\includegraphics[width=5.5cm]{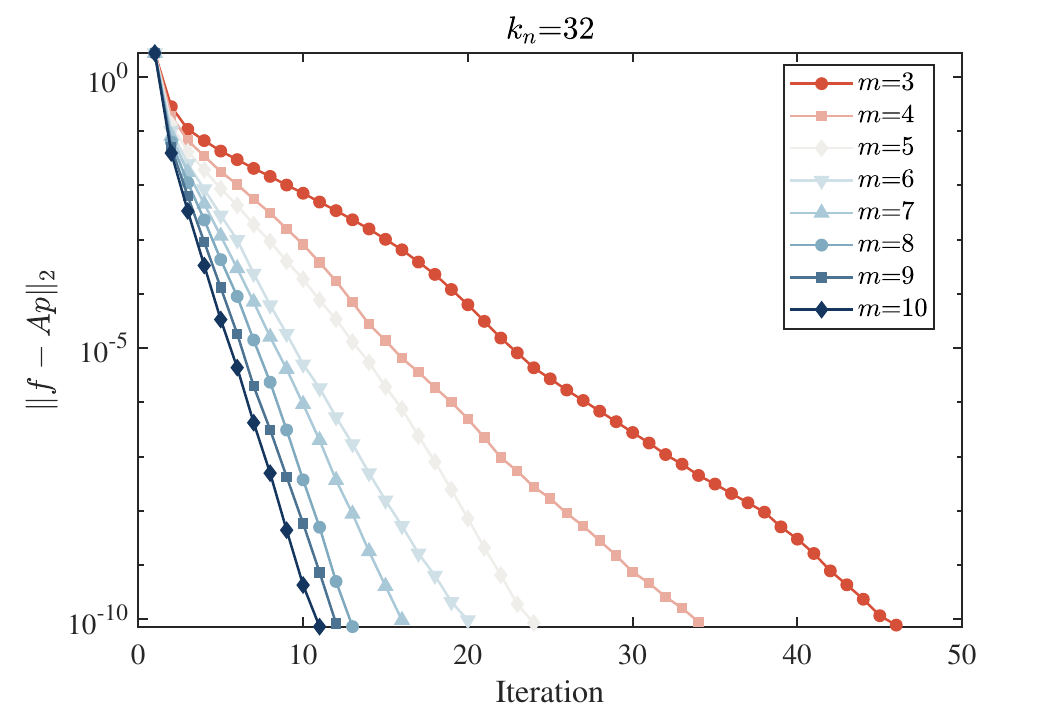}}

	\subfloat[\centering]{\includegraphics[width=5.5cm]{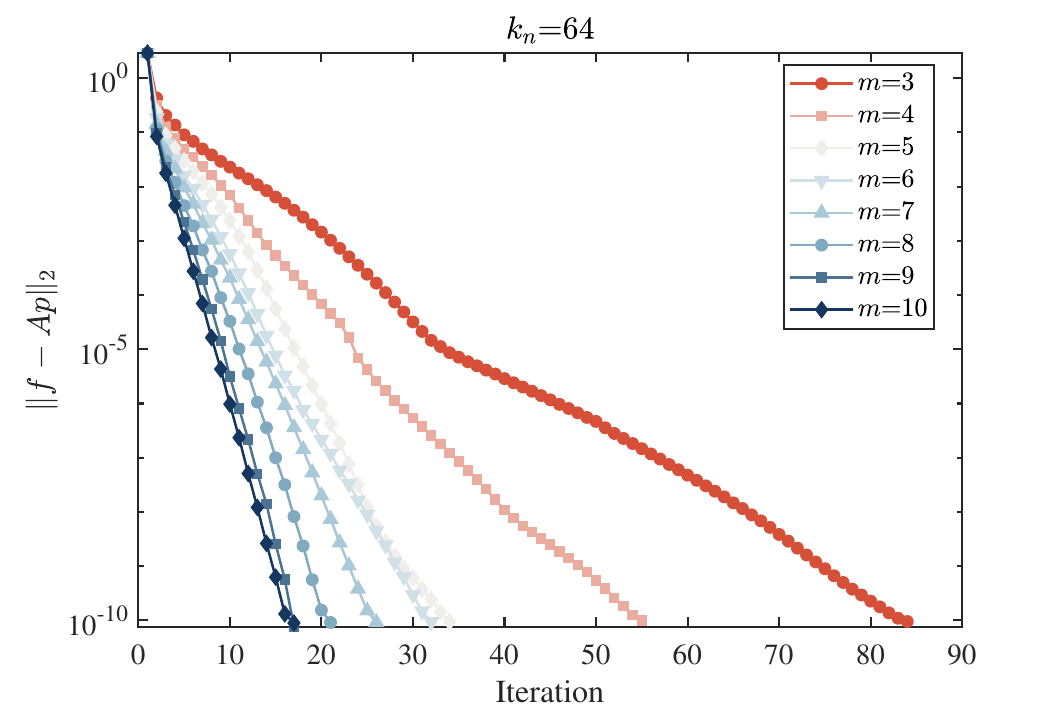}}
	\subfloat[\centering]{\includegraphics[width=5.5cm]{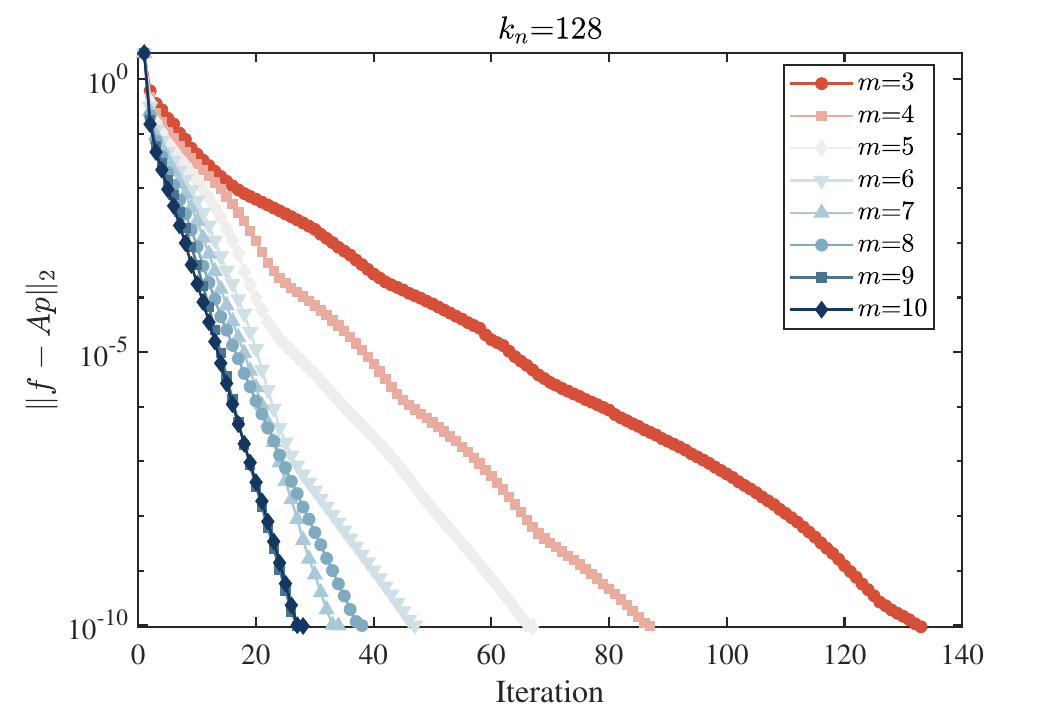}}
	\caption{The iterative residual reduction in GMRES iteration with different size of the equation $k_n$ and different dimension of subspace of GMRES $m$.}
	\label{fig:iter_res}
\end{figure}

There exists a relationship between grid resolution and the optimal GMRES subspace dimension for efficient convergence. To assess the impact of grid number and subspace dimension on the algorithm, we evaluated the effects of different grid numbers and subspace dimensions on the convergence of GMRES iterations, as shown in Figure \ref{fig:iter_res}.

\begin{figure}[H]
    \centering
    \subfloat[\centering]{\includegraphics[width=6.0cm]{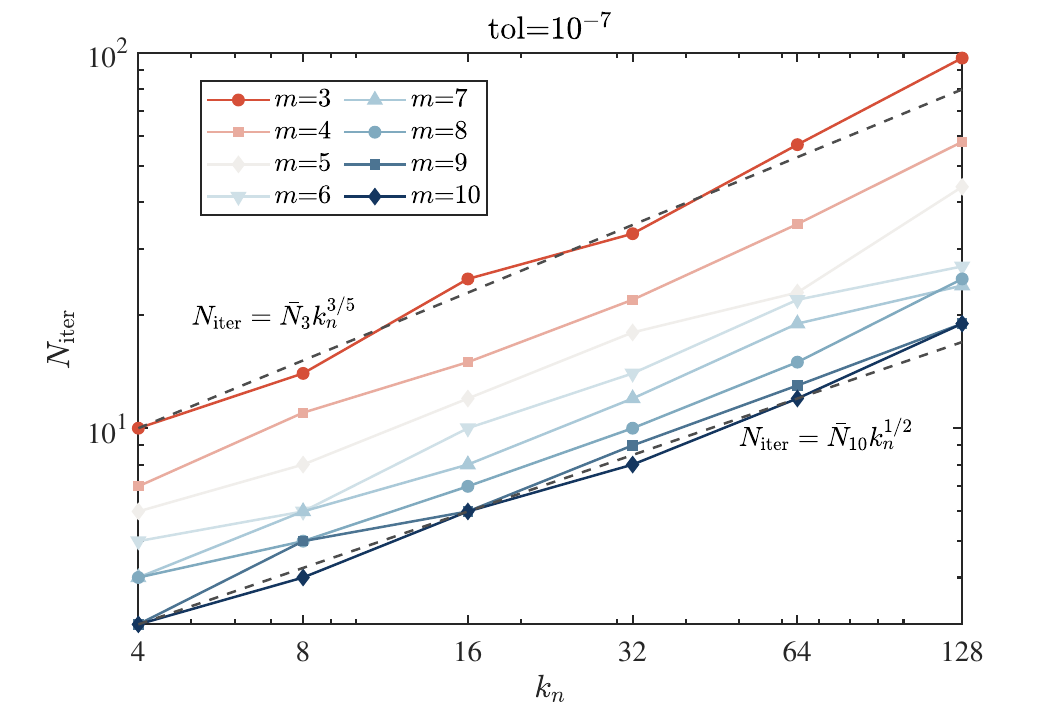}}
    \subfloat[\centering]{\includegraphics[width=6.0cm]{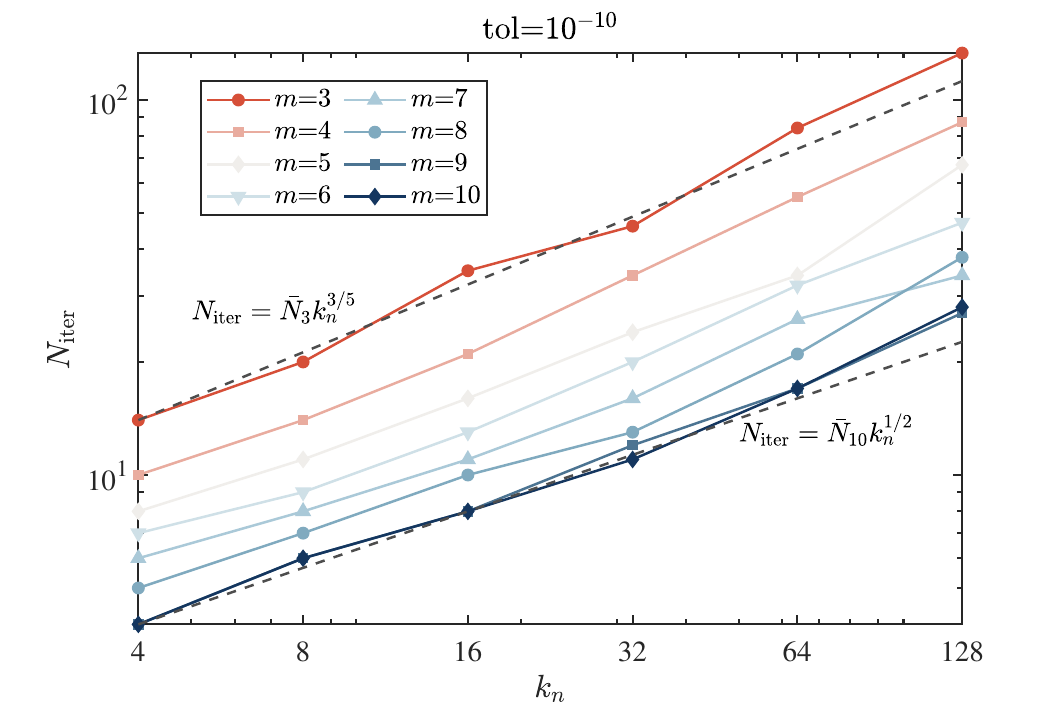}}
    \caption{The iteration number for convergence under different residual tolerances.}
    \label{fig:cross_iter_count}
\end{figure}

As shown in Figure \ref{fig:iter_res}, the iteration count increases with $k_n$ (matrix size), a trend similar to that observed in other GMRES-like algorithms. It is well known that the FFT on simple geometries offers a direct solution without iterations; therefore, its computational cost remains consistent regardless of the matrix size. Given that the total computational cost of an iterative algorithm is determined by both the number of iterations and the computational cost per iteration, we next evaluate the impact of matrix size on the iteration count, and the results are presented in Figure \ref{fig:cross_iter_count}.

Figure \ref{fig:cross_iter_count} illustrates the iteration counts for two representative residual convergence criteria, $10^{-7}$ and $10^{-10}$. Two power-law fitting curves, depicted as dashed lines in Figure \ref{fig:cross_iter_count}, are introduced as references of computational complexity order.  Note that $\bar{N}_3$ and $\bar{N}_{10}$ are fitting constants independent of $k_n$.  As can be seen, regardless of the subspace dimension $m$, the iteration count scales approximately with $k_n^{1/2}$. Since the computational complexity of each iteration is $O\left(k_n^2\log\left(k_n\right)\right)$, the estimated total computational complexity is $O\left(k_n^{5/2}\log\left(k_n\right)\right)$. This represents a significant improvement compared to the complexity of traditional methods, which is typically higer than $O\left(k_n^3\right)$. The actual computation time of the algorithm is further assessed below to demonstrate its computational efficiency.\par

\subsection{Computational efficiency}
\label{ssec:effi}

\begin{figure}[H]
    \centering
    \subfloat[\centering]{\includegraphics[width=6.0cm]{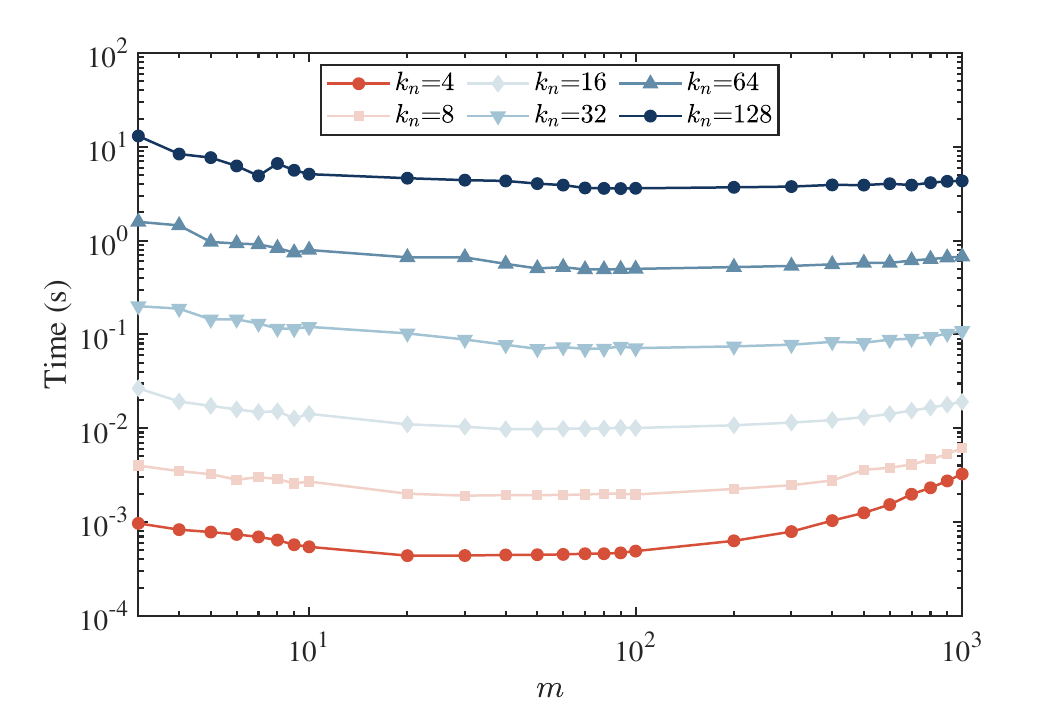}}
    \subfloat[\centering]{\includegraphics[width=6.0cm]{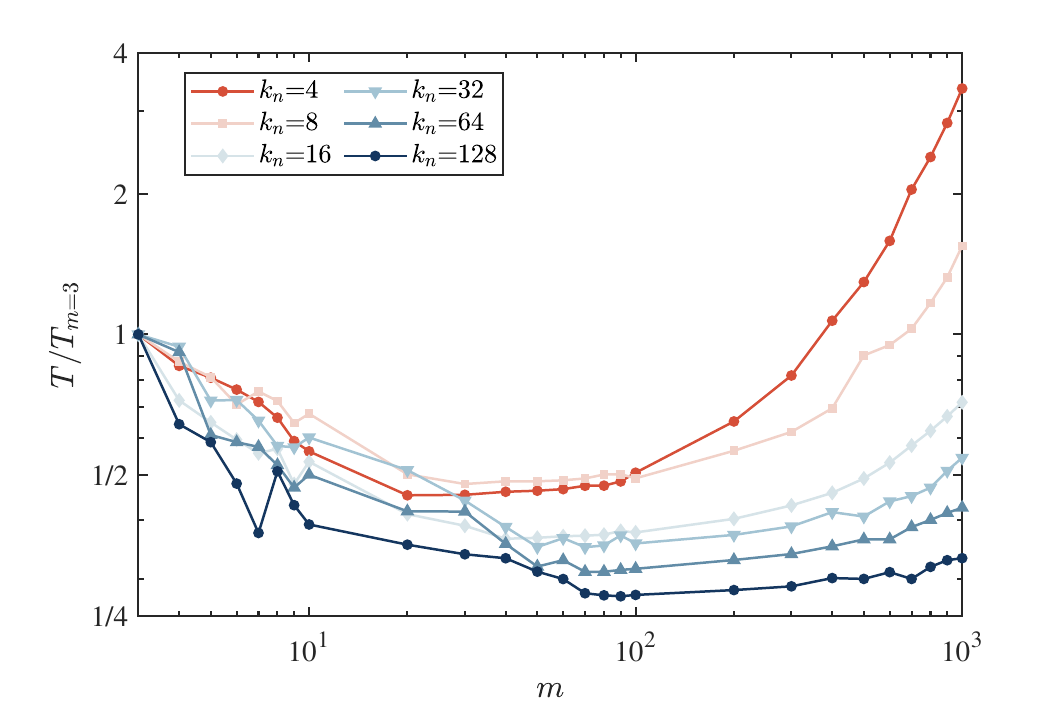}}
    \caption{Computational time with different subspace dimension $m$ and equation size $k_n$.}
    \label{fig:cross_m_time}
\end{figure}

The results in Figure \ref{fig:iter_res} indicate that the subspace dimension $m$ has a significant impact on the algorithm, specifically in three aspects: 1. When the subspace dimension $m$ is small, the number of iterations increases, inevitably slowing down the computational efficiency. 2. When the subspace dimension is large, the time for each GMRES iteration increases. An excessively large subspace dimension $m$ can even lead to an increase in the computational cost of solving the coefficient equation within the subspace itself. Therefore, to evaluate the impact of the subspace dimension $m$ on computational efficiency, we assessed the influence of the subspace dimension $m$ on the GMRES computation time. The results are shown in Figure \ref{fig:cross_m_time}.

In Figure \ref{fig:cross_m_time}(a), the impact of the subspace dimension $m$ on the total GMRES computation time is presented for six different equation sizes, denoted by $k_n$. It can be observed that for each matrix size, the influence of the subspace dimension $m$ on the total GMRES time aligns with theoretical expectations; specifically, both excessively small and excessively large subspace dimensions $m$ lead to increased computation time. To more clearly illustrate the effect of the subspace dimension $m$ on computation time, Figure \ref{fig:cross_m_time}(b) provides a comparison of computation times at various subspace dimensions $m$ against a baseline. This baseline is defined as the GMRES computation time for each equation size $k_n$ when $m=3$. It is clearly evident that when the subspace dimension $m$ is approximately less than 80, the computation time decreases as $m$ increases. Conversely, when the subspace dimension $m$ exceeds 100, the computation time rapidly increases with increasing $m$. Consequently, for all parallel efficiency tests, $m=80$ was consistently chosen as the optimal subspace dimension.\par

\begin{figure}[H]
    \centering
    \subfloat[\centering]{\includegraphics[width=6.0cm]{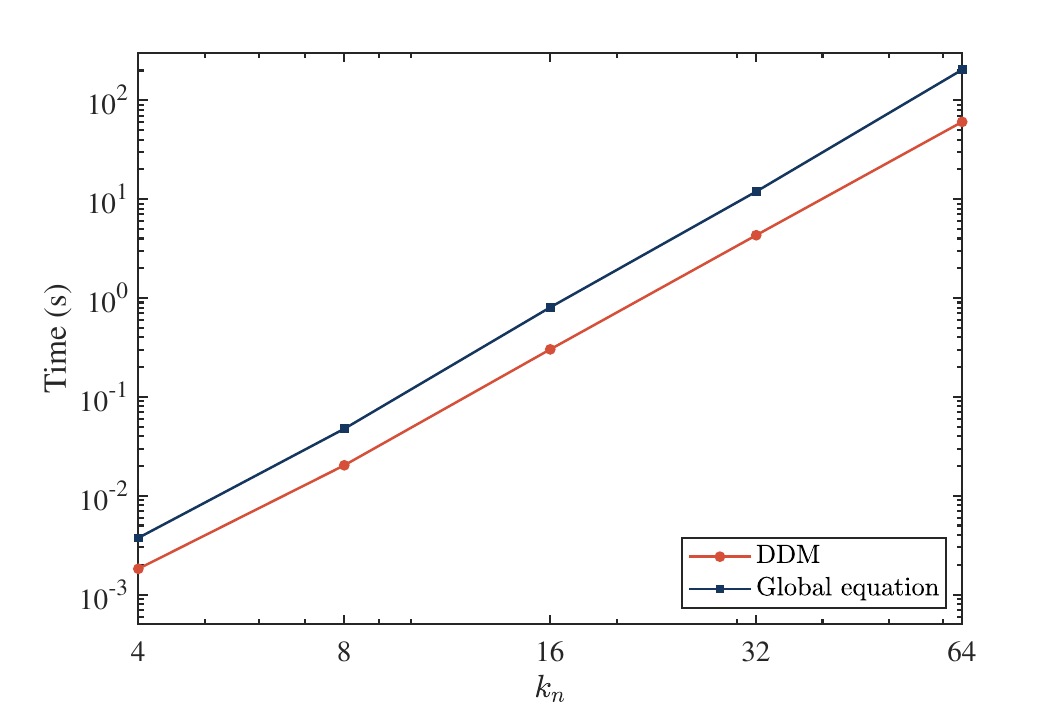}}
    \subfloat[\centering]{\includegraphics[width=6.0cm]{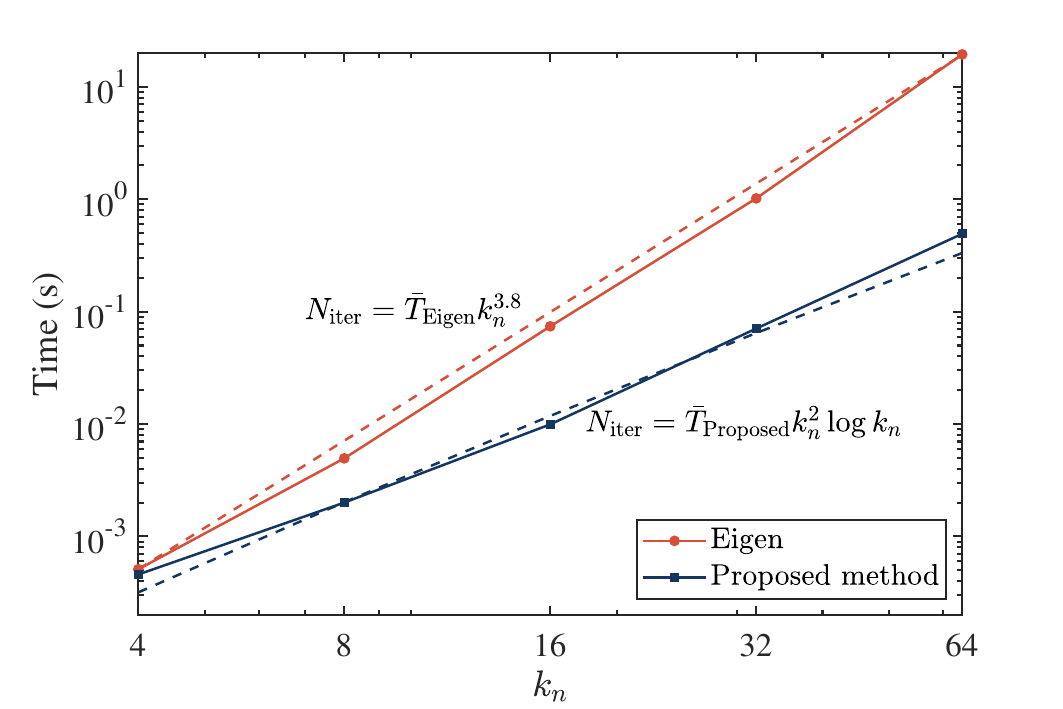}}
    \caption{Computational time comparison for the cross-shaped domain. (a) Total solution time for the DDM solver vs. a direct global solve, both based on Eigen. (b) Solution time for a single sub-domain using the FFT-based GMRES solver (Prposed method) vs. Eigen solver.}
    \label{fig:cross_time_comp}
\end{figure}

To validate the performance of the proposed algorithm, we conduct two numerical experiments: First, we compare the total solution time of our DDM-based solver against a direct solve of the global system for the cross-shaped domain. Since the complex geometry precludes a direct application of FFT-based methods to the global problem, we employ a sparse GMRES solver from the Eigen C++ library \cite{eigenweb} as a benchmark. This comparison, presented in Figure \ref{fig:cross_time_comp}(a), is designed to demonstrate the overall effectiveness of the DDM strategy. Second, to specifically evaluate the efficiency of our fast sub-domain solver, we compare its performance against Eigen's sparse GMRES solver on a single rectangular sub-domain (taken as $\Omega_2$). This test isolates the performance of the GMRES solver component from the overhead of the DDM framework. The results of this direct solver comparison are shown in Figure \ref{fig:cross_time_comp}(b).

In this comparison, all solvers were evaluated under identical operating environments, encompassing both hardware and parallel computing settings. Additionally, the Eigen solver employed sparse matrix storage and computation techniques. As shown in Figure \ref{fig:cross_time_comp}(a), for any given equation size, the computational efficiency of the domain decomposition-based algorithm is several times greater than that of directly solving the global system of equations. This demonstrates that the adoption of a domain decomposition strategy significantly enhances computational efficiency compared to directly solving the global system using GMRES.\par

Furthermore, Figure \ref{fig:cross_time_comp}(b) illustrates that the proposed algorithm exhibits a notable computational efficiency advantage over Eigen for large-scale problems. This superiority stems from differences in computational complexity, which aligns with our prior analysis. To quantify this performance advantage, the solution times plotted in Figure \ref{fig:cross_time_comp}(b) are fitted to curves (dashed lines) that represent the observed computational scaling, where $\bar{T}_{\rm Eigen}$ and $\bar{T}_{\rm Proposed}$ are fitting constants. A comparison reveals that the computation time for traditional GMRES and similar solvers, when evaluated with respect to $k_n$, scales approximately as $k_n^{3.8}$. In contrast, the solution time for the algorithm proposed in this study scales approximately as $k_n^{2}\log k_n$. Both trends closely match their respective theoretical estimations. While $k_n=64$, which corresponds to a total of approximately $1.4\times 10^5$ grid points, the proposed algorithm achieves a 40-times speedup compared to the Eigen-based solver.\par

Furthermore, due to fundamental limitations inherent in the computational principles of Krylov subspace methods, optimizing complexity cannot be achieved merely by replacing preconditioners without altering the underlying per-iteration algorithm or other core algorithmic components. Consequently, the GMRES algorithm equipped with an FFT-based preconditioner, as proposed in this study, offers substantial advantages for computations involving large matrices.

\section{Conclusion}
\label{sec:conclusion}
In this study, we propose a GMRES solution algorithm employing an FFT-based preconditioner for solving the Poisson equation in complex geometries. This algorithm achieves an $O(N \log N)$ complexity for Poisson equation solutions within complex domains, demonstrating significant computational advantages for problems involving large mesh sizes. The specific contributions of this work are as follows:

\begin{enumerate}
    \item A matrix equation-based FFT-type algorithm framework was constructed. Unlike traditional FFT-based Poisson equation solvers that often employ spectral space expansion strategies, this framework presents an accelerated FFT-based solution algorithm leveraging eigenvalue transformations.
    \item Building upon the matrix equation-based FFT-type algorithm framework, the integration of FFT algorithms with DDMs was first realized. This extends the applicability of FFT-type Poisson equation solvers to complex geometries formed by concatenating multiple regions.
    \item A GMRES algorithm with an FFT-based preconditioner was proposed. Compared to traditional preconditioners, this significantly reduces the number of GMRES iterations for sub-equations. Furthermore, an assessment of the GMRES iteration count indicates that the complexity of the FFT-based Poisson equation solver is only increased to $O(k_n^{5/2}\log k_n)$, which still presents a significant advantage over traditional algorithms with $O(k_n^3)$ computational complexity, such as the standard GMRES algorithm.
    \item An in-house C++ parallel solver was developed based on FFTW. Its computational efficiency was compared with the Eigen solver on a parallel platform. The results demonstrate that the developed algorithm exhibits an $O(k_n^2\log k_n)$ computational complexity, which is significantly superior to the $O(k_n^4)$ computational complexity of traditional GMRES algorithms. For GMRES solutions of matrices with a scale of $2^{15}\times2^{15}$, the algorithm proposed in this study achieves a computational efficiency improvement of more than 20 times.
\end{enumerate}

\section*{Acknowledgement}
This work was supported by Guangdong Science and Technology Fund, Grant 2022B1515120009. No external datasets, either public or private, were utilized in this study; all computational results were generated solely through the methods described. The in-house solver developed and mentioned in this research has been open-sourced and is publicly available at \url{https://github.com/ffskibkwi/Focalors_Poisson}. All computational test results presented in this paper are fully reproducible using the algorithms detailed herein.\par

During the preparation of this work the authors used Gemini in order to improve the English language, correct grammar and punctuation, and enhance the clarity and conciseness of the text, without generating any scientific content or ideas. After using this tool, the authors reviewed and edited the content as needed and take full responsibility for the content of the published article.



\bibliography{mybibfile}

\end{document}